# ON THE SUPER REPLICATION PRICE OF UNBOUNDED CLAIMS

By Sara Biagini and Marco Frittelli

*Università degli Studi di Perugia and Università degli Studi di Firenze*

In an incomplete market the price of a claim $f$ in general cannot be uniquely identified by no arbitrage arguments. However, the "classical" super replication price is a sensible indicator of the (maximum selling) value of the claim. When $f$ satisfies certain pointwise conditions (e.g., $f$ is bounded from below), the super replication price is equal to $\sup_Q E_Q[f]$, where $Q$ varies on the whole set of pricing measures. Unfortunately, this price is often too high: a typical situation is here discussed in the examples.

We thus define the less expensive *weak super replication price* and we relax the requirements on $f$ by asking just for "enough" integrability conditions.

By building up a proper duality theory, we show its economic meaning and its relation with the investor's preferences. Indeed, it turns out that the weak super replication price of $f$ coincides with $\sup_{Q \in M_\Phi} E_Q[f]$, where $M_\Phi$ is the class of pricing measures with finite generalized entropy (i.e., $E[\Phi(\frac{dQ}{dP})] < \infty$) and where $\Phi$ is the convex conjugate of the utility function of the investor.

**1. Introduction.** We investigate the super replication price of contingent claims in incomplete markets where gains from trading may take any real value. For claims $f$ which are bounded from below, the classical super replication price is equal to

$$\sup_{Q \in M_1} E_Q[f], \tag{1}$$

where $M_1$ is the set of all pricing measures. For claims which are unbounded from below, however, the above supremum may be strictly lower than the super replication price.

One of the main results of the paper is a representation of the supremum (1) for unbounded claims in terms of a "weak super replication price" $\hat{f}_\Phi$, which allows variables from a slightly wider class than the usual one of terminal values from admissible integrands. This natural class $C_\Phi$ (see [15]) was first explicitly introduced by Frittelli (see [8, 9]). The class $C_\Phi$ depends on a convex function $\Phi:(0,+\infty) \to \mathbb{R}$ which normally (see Remark 7) represents the conjugate function of a utility function $u$. We will assume that $\Phi$ satisfies a growth condition that is shown to be equivalent to the condition of reasonable asymptotic elasticity of $u$ in the sense of Schachermayer [19].

We denote by $M_\Phi \triangleq \{Q \in M_1 : E[\Phi(\frac{dQ}{dP})] < \infty\}$ the set of pricing measures with finite generalized entropy. The actual result obtained (see Theorem 5) is that if $\Phi(0) < \infty$ and there exists an equivalent pricing measure with finite generalized entropy, then for claims $f$ (for which the LHS make sense, but which may be unbounded from below) we have

$$(2) \qquad \sup_{Q \in M_\Phi} E_Q[f] = \inf\{x \in \mathbb{R} | f - x \in C_\Phi\} \triangleq \hat{f}_\Phi.$$

The representation of (1) is then a corollary, setting $\Phi = id$.

We provide an example of an unbounded claim where the weak super replication price $\hat{f}_{id}$ is strictly less than the classical super replication price $\hat{f}$.

The paper is based on the appropriate selection of the spaces for which the following duality holds true: if $\Phi(0) < \infty$ (and there exists an equivalent pricing measure in $M_\Phi$), then the cones $C_\Phi$ and $co(M_\Phi)$ are polar to one another.

However, if $\Phi(0)$ is infinite, then $co(M_\Phi) \subseteq (C_\Phi)^0$ with possibly strict inclusion. We give an example where indeed the inclusion is strict and $co(M_\Phi)$ is not closed.

Finally, we develop a comparison between the duality relation obtained by Delbaen and Schachermayer [5] and ours when $\Phi = id$. It turns out that the super replication price $\hat{f}_w$ of the claim $f$, as defined in [5], depends explicitly on an unbounded weight function $w$, which represents the maximum loss the investor is willing to face. Instead, our weak super replication price $\hat{f}_{id}$ is equal for all the agents in the given market.

If one is interested in taking into account the investor's attitude toward risk, we suggest $\hat{f}_\Phi$ as a suitable super replication price, since it has the advantage of being explicitly linked to the utility function.

The paper is organized as follows.

Section 1 has three sections: the first contains the general setup and the precise formulations of our results; in the second we explain how the preferences of the investors are taken into consideration and the relations between $u$ and $\Phi$; the third is devoted to two basic examples in which classical duality fails.

In Section 2 we give an abstract duality relation, which is used in the proofs of the main results, and we also provide a new proof of the representation of the super replication price for bounded-from-below claims.

In Section 3 we build up a proper dual system, so that we obtain the polarity between $C_\Phi$ and $co(M_\Phi)$ and we prove (2).

We end with Section 4, which contains the comparison between $\hat{f}_{id}$ and $\hat{f}_w$.

1.1. *The model and the results.* Our starting point is the general semimartingale model of a financial market as defined by Delbaen and Schachermayer [5].

Let $(\Omega, \mathcal{F}, (\mathcal{F}_t)_{t \in [0,T]}, P)$ be a filtered probability space, where we assume that the filtration satisfies the usual assumptions of right continuity and completeness, and let $\mathbb{P}$ be the class of probability measures equivalent to $P$.

The $\mathbb{R}^d$-valued càdlàg semimartingale $X = (X_t)_{t \in [0,T]}$ represents the (discounted) price process of $d$ tradeable assets.

An $\mathbb{R}^d$-valued predictable process $H = (H_t)_{t \in [0,T]}$ is called an *admissible* trading strategy if $H$ is $X$-integrable and there exists a constant $c \in \mathbb{R}$ such that, for all $t \in [0,T]$, $\int_0^t H_s \cdot dX_s \geq -c$, $P$-a.s. The financial interpretation of $c$ is a finite credit line which the investor must respect in his or her trading. This bounded-from-below restriction on the stochastic integral traces back to the work of Harrison and Pliska [13] and it is now a standard assumption in the literature (see [4]).

We denote by $L^0$ [resp. $L^\infty$, $L^1(P)$] the space of $P$-a.s. finite (resp. $P$-essentially bounded, $P$-integrable) random variables on $(\Omega, \mathcal{F})$, with $L^\infty_+$ (resp. $L^1_+$) the cone of $P$-a.s. nonnegative random variables in $L^\infty$ (resp. $L^1$), with $L^{bb}$ the cone of essentially bounded from below random variables, with $\overline{C}^P$ the closure of a set $C \subseteq L^1(P)$ in the $L^1(P)$ norm topology. Define

$$K \triangleq \left\{ \int_0^T H_s \cdot dX_s \middle| H \text{ is admissible} \right\} \subseteq L^{bb},$$

$$C \triangleq (K - L^0_+) \cap L^\infty.$$

$K$ is the cone of all claims that are replicable, at zero initial cost, via admissible trading strategies. The set

$$(K - L^0_+) = \{ f \in L^0 : \exists g \in K \text{ s.t. } g \geq f \ P\text{-a.s.} \}$$

is the cone of all claims in $L^0$ that can be dominated by a replicable claim, hence is the cone of super-replicable claims. Consequently $C \triangleq (K - L^0_+) \cap L^\infty$ is the cone of bounded super-replicable claims. In Section 3 we will consider the closure $\overline{C}$ of $C$ under a particular topology: then $\overline{C}$ is the cone of claims that can be "approximated" by bounded super-replicable claims.

Define

(3) $\quad M_1 \triangleq \{Q \ll P : K \subseteq L^1(Q) \text{ and } E_Q[g] \leq 0 \text{ for all } g \in K\},$

(4) $\quad M \triangleq \{z \in L^1(P) : E[zg] \leq 0 \ \forall g \in C\} \subseteq L^1_+(P).$

The elements in $M_1$ are called *separating probability measures*. We will often identify probability measures $Q$, absolutely continuous with respect to $P$, with their Radon–Nikodym derivatives $\frac{dQ}{dP} \in L^1(P)$. Note that (see [2], Lemma 1.1 for details)

(5)
$$M_1 = \{Q \ll P : E_Q[g] \leq 0 \ \forall g \in C\}$$
$$= \{z \in M | E[z] = 1\}$$

and that if $X$ is bounded (resp. locally bounded), then

$$M_1 = \{Q \ll P : X \text{ is a } (Q, (\mathcal{F}_t)_{t \in [0,T]}) \text{ martingale (resp. local martingale)}\},$$

that is, $M_1$ is the set of $P$-absolutely continuous *martingale* (resp. *local martingale*) measures. In general, for possibly unbounded $X$, $M_1$ is the set of $P$-absolutely continuous probabilities such that the admissible stochastic integrals are supermartingales. What is more (see [5], Proposition 4.7) if $M_1 \cap \mathbb{P} \neq \varnothing$, then the set $M_\sigma$ of absolutely continuous $\sigma$-*martingale* probabilities is not empty and $M_\sigma$ is dense in $M_1$ for the total variation topology.

The main topic of this paper is the analysis of the *super replication price* $\hat{f}$ of a claim $f \in L^0$, defined by

$$\hat{f} \triangleq \inf\{x \in \mathbb{R} | \exists g \in K \text{ s.t. } x + g \geq f \ P\text{-a.s.}\}$$
$$= \inf\{x \in \mathbb{R} | f - x \in (K - L^0_+)\}.$$

This subject was originally studied by El Karoui and Quenez [7]; see also Karatzas [15] and the references cited there. We will mainly deal with the results on this subject provided by Delbaen and Schachermayer **(year?)**. If $f \in L^1(Q)$ for all $Q \in M_1$, then

(6)
$$\hat{f} = \inf\left\{x \in \mathbb{R} \Big| f - x \in (K - L^0_+) \bigcap_{Q \in M_1} L^1(Q)\right\}$$
$$= \inf\left\{x \in \mathbb{R} \Big| f - x \in \bigcap_{Q \in M_1} (K - L^1_+(Q))\right\}$$

since, for all $Q \in M_1$, $(K - L^0_+) \cap L^1(Q) = (K - L^1_+(Q))$.

If $f \in L^{bb}$, then

$$\hat{f} = \inf\{x \in \mathbb{R} | f - x \in (K - L^0_+) \cap L^{bb}\} = \inf\{x \in \mathbb{R} | f - x \in C_{bb}\},$$

where
$$C_{bb} \triangleq (K - L_+^0) \cap L^{bb}.$$

It is easy to see that $\hat{f}$ dominates $\sup_{Q \in M_1} E_Q[f]$.

PROPOSITION 1. *If $M_1 \neq \varnothing$ and if either $f \in \bigcap_{Q \in M_1} L^1(Q)$ or $f \in L^{bb}$, then*

(7) $$\sup_{Q \in M_1} E_Q[f] \leq \hat{f}.$$

PROOF. For all $x \in \mathbb{R}$ such that $f - x \in (K - L_+^0)$ we have $0 \geq \sup_{Q \in M_1} E_Q[f - x] = \sup_{Q \in M_1} E_Q[f] - x$. □

REMARK 2. If $N$ is a convex set of probability measures absolutely continuous with respect to $P$ and if $N \cap \mathbb{P} \neq \varnothing$, then it is easy to show that if $f \in \bigcap_{Q \in N} L^1(Q)$ or if $f \in L^{bb}$, then

(8) $$\sup_{Q \in N} E_Q[f] = \sup_{Q \in N \cap \mathbb{P}} E_Q[f].$$

In fact, let $Q_0 \in N$ and $Q_1 \in N \cap \mathbb{P}$: take the convex combinations $Q^x = (1-x)Q_0 + xQ_1$, $x \in [0,1]$. If $x \to 0$, then $\frac{dQ^x}{dP} \to \frac{dQ_0}{dP}$ in $L^1(P)$ and also $P$-almost surely. In case $f \in L^{bb}$, equality (8) is a simple consequence of Fatou's lemma. In case $f \in \bigcap_{Q \in N} L^1(Q)$, we have $|f|\frac{dQ^x}{dP} \leq |f|(\frac{dQ_0}{dP} + \frac{dQ_1}{dP})$ and so the dominated convergence theorem can be applied. Therefore, in what follows (Theorem 3, Corollary 4, Theorem 5 and Proposition 6) it will be equivalent to take the supremum over the sets $M_1$ ($M_\Phi$) or over $M_1 \cap \mathbb{P}$ ($M_\Phi \cap \mathbb{P}$).

Delbaen and Schachermayer proved ([5], Theorem 5.10) that in (7) equality holds if $f$ is bounded from below:

THEOREM 3. *If $M_1 \cap \mathbb{P} \neq \varnothing$ and if $f \in L^{bb}$, then*

(9) $$\hat{f} = \sup_{Q \in M_1} E_Q[f].$$

A new proof of this result is given in Section 2.1.

If $f \in \bigcap_{Q \in M_1} L^1(Q)$, (9) does not hold true anymore, when $\hat{f}$ is given in (6). To obtain a correct dual formula, we must replace in (6) the set $\bigcap_{Q \in M_1}(K - L_+^1(Q))$ with $\bigcap_{Q \in M_1} \overline{K - L_+^1(Q)}^Q \triangleq C_{id}$, that is, with the closure of $C$ under an appropriate topology (see Theorem 17). As a consequence of Theorem 5 below, with $\Phi = id$, we deduce the following.

COROLLARY 4. *If $M_1 \cap \mathbb{P} \neq \varnothing$ and if $f \in \bigcap_{Q \in M_1} L^1(Q)$, then*

$$\text{(10)} \quad \hat{f}_{id} \triangleq \inf\left\{x \in \mathbb{R} \,\Big|\, f - x \in \bigcap_{Q \in M_1} \overline{K - L^1_+(Q)}^Q\right\} = \sup_{Q \in M_1} E_Q[f].$$

We shall call $\hat{f}_{id}$ the *weak super replication price* of $f$. In Example 8 of Section 1.3 we show that it is possible that $\hat{f}_{id} < \hat{f}$.

The introduction of the convex function $\Phi$ will allow us to present our results in a more general framework and to link the interpretation of the weak super replication price with the preferences of an investor represented by his or her utility function. This analysis is provided in Section 1.2.

Throughout the paper we make the following assumption.

ASSUMPTION. The function $\Phi:(0,+\infty) \to \mathbb{R}$ is convex and satisfies the following growth condition:

$G(\Phi): \forall [\lambda_0, \lambda_1] \subseteq (0,+\infty)$ there exist $\alpha > 0, \beta > 0$ such that
$$\Phi^+(\lambda y) \leq \alpha \Phi^+(y) + \beta(y+1) \;\forall\, y > 0, \;\forall\, \lambda \in [\lambda_0, \lambda_1].$$

For a detailed discussion of this condition and its relation with the condition, introduced by Schachermayer [19], of *reasonable asymptotic elasticity* of the utility function we defer to [10]. Set $\Phi(0) = \lim_{y \downarrow 0} \Phi(y)$ and define:

$$M_\Phi \triangleq \left\{Q \in M_1 : \Phi\left(\frac{dQ}{dP}\right) \in L^1(P)\right\}.$$

In Example 8, where $\Phi$ is the identity function *id* and so $M_\Phi = M_1$, we will show that if $f \in \bigcap_{Q \in M_\Phi} L^1(Q)$, then it may happen that

$$\inf\left\{x \in \mathbb{R} \,\Big|\, f - x \in \bigcap_{Q \in M_\Phi} (K - L^1_+(Q))\right\} > \sup_{Q \in M_\Phi} E_Q[f].$$

The examples in Section 1.3 and the next theorem, proved in Section 3, are the main contributions of the paper. Our aim is exactly that of providing the correct interpretation and the dual representation of $\sup_{Q \in M_\Phi} E_Q[f]$, even when it is strictly less than $\hat{f}$.

THEOREM 5. *If $\Phi(0) < +\infty$, $M_\Phi \cap \mathbb{P} \neq \varnothing$ and $f \in \bigcap_{Q \in M_\Phi} L^1(Q)$, then*

$$\text{(11)} \quad \hat{f}_\Phi \triangleq \inf\left\{x \in \mathbb{R} \,\Big|\, f - x \in \bigcap_{Q \in M_\Phi} \overline{K - L^1_+(Q)}^Q\right\} = \sup_{Q \in M_\Phi} E_Q[f] \leq \hat{f}.$$

As already mentioned, in Theorem 17 we will show that $\bigcap_{Q \in M_\Phi} \overline{K - L^1_+(Q)}^Q = \overline{C} = C_\Phi$, where $\overline{C}$ is the closure of $C$ under an appropriate topology.

As a consequence of Theorem 1.1 of Kabanov and Stricker [14] we also have

PROPOSITION 6. *If $M_\Phi \cap \mathbb{P} \neq \varnothing$ and $f \in L^{bb}$, then*
$$\hat{f} = \sup_{Q \in M_1} E_Q[f] = \sup_{Q \in M_\Phi} E_Q[f] = \hat{f}_\Phi.$$

PROOF. By definition, if $f \in L^{bb}$, then $\hat{f}_\Phi \leq \hat{f}$. As in the proof of Proposition 1 we also get $\sup_{Q \in M_\Phi} E_Q[f] \leq \hat{f}_\Phi$. The growth condition $G(\Phi)$ is weaker than the condition used in Corollary 1.4 of [14], since $G(\Phi)$ does not require that $\Phi(0) < +\infty$. Nevertheless, it can be shown, as in the proof of Corollary 1.4 of [14], that the condition $G(\Phi)$ and Theorem 1.1 of [14] imply

(12) $$\sup_{Q \in M_\Phi} E_Q[f] = \sup_{Q \in M_1} E_Q[f] \qquad \text{if } f \in L^{bb}.$$

Hence, from (9), we get $\hat{f} = \sup_{Q \in M_1} E_Q[f] = \sup_{Q \in M_\Phi} E_Q[f] \leq \hat{f}_\Phi \leq \hat{f}$. □

In Example 9 we will show that the equality $\hat{f}_\Phi = \hat{f}$ may not be true for claims that are not bounded from below.

1.2. *Taking preferences into account.* In incomplete markets, it may be useful to take into account the preferences of the investor. This naturally leads to the specification of a utility function $u$, which we assume to be strictly concave, increasing and finite valued on the whole $\mathbb{R}$. The related standard utility maximization problem

$$\sup_{g \in K} E[u(x + g)], \qquad x \in \mathbb{R},$$

in general does not admit an optimal solution in $K$ (see [19]). In the duality theory approach to this problem a crucial role is played by the convex conjugate of $u$, which we denote by $\Phi$:

$$\Phi(y) \triangleq \sup_{x \in \mathbb{R}} \{u(x) - xy\}, \qquad y > 0.$$

Note that the condition $\Phi(0) < +\infty$ assumed in Theorem 5 is equivalent to the requirement that the utility function is bounded from above.

REMARK 7. The function $\Phi = id$ is the convex conjugate of the function $u : \mathbb{R} \to \mathbb{R} \cup \{-\infty\}$ defined by

$$u(x) = \begin{cases} 0, & \text{if } x = -1, \\ -\infty, & \text{otherwise,} \end{cases}$$

which is not increasing on $\mathbb{R}$. In this case $\Phi$ cannot be interpreted as the conjugate of a "utility" function.

It was first shown in [2] that if

$$\sup_{g \in K} E[u(x+g)] < u(+\infty),$$

then the fundamental duality relation

$$\sup_{g \in K} E[u(x+g)] = \min_{Q \in M_\Phi} \min_{\lambda > 0} \lambda x + E\left[\Phi\left(\lambda \frac{dQ}{dP}\right)\right]$$

holds true, without any further assumption on the utility function. For what concerns economic considerations, Frittelli [9] suggested a clear financial interpretation for the class $M_\Phi$ of those separating measures having finite *generalized entropy*. In fact, fix $Q \in M_1$ and consider the problem

$$U_Q(x) \triangleq \sup\{E[u(x+g)] | g \in L^1(Q),\ E_Q[g] \leq 0,\ u^-(x+g) \in L^1(P)\}.$$

This is precisely the utility maximization problem we would face if we selected $Q$ as pricing measure. When $G(\Phi)$ is satisfied, then (see [9], Proposition 4) $Q$ belongs to $M_\Phi$ if and only if

$$U_Q(x) < u(+\infty) \qquad \text{for all } x \in \mathbb{R}.$$

More explicitly this means that pricing by $Q \in M_\Phi$ guarantees that the investor cannot reach his or her maximum possible utility, $u(+\infty)$, starting with an arbitrarily low initial endowment $x$. Therefore it makes sense to work with $M_\Phi$, as the class of pricing measures which makes the model free of this types of *utility based arbitrage opportunities*.

1.3. *Examples.* In Example 8 we show that $\hat{f}_{id} < \hat{f}$ and in Example 9 we show a case where $\hat{f}_\Phi < \hat{f}$, when $\Phi$ is not the identity function.

EXAMPLE 8. We denote by $I_n$ the interval $(\frac{1}{2^n}, \frac{1}{2^{n-1}}]$ and by $J_n^1$ and by $J_n^2$ its two halves $(\frac{1}{2^n}, \frac{3}{2^{n+1}}]$ and $(\frac{3}{2^{n+1}}, \frac{1}{2^{n-1}}]$, respectively.

We consider the following one-period model: $(\Omega, (\mathcal{F}_0, \mathcal{F}_1), P)$, where $\Omega$ is the interval $(0, 1]$, $\mathcal{F}_0 = \sigma\{I_n | n \in \mathbb{N}_0\}$, $\mathcal{F}_1 = \sigma\{J_n^i | i = 1, 2 \text{ and } n \in \mathbb{N}_0\}$ and $P$ is the restriction of the Lebesgue measure to $\mathcal{F}_1$. The process $X$ is given by $X(0) = 0$ and

$$X(1) = \begin{cases} n, & \text{on } J_n^1, \\ -n^2, & \text{on } J_n^2. \end{cases}$$

The set $K^0$ will be the set of all stochastic integrals with respect to predictable processes, with no admissibility restrictions. Here this set is simply $\{\alpha X(1) | \alpha\ \mathcal{F}_0\text{-measurable}\}$ and $\alpha$ is identified by the sequence $(\alpha_n)_{n \geq 1}$ of its values on the intervals $I_n$. The structure of elements in $K$ can now be easily described. By fixing a credit level $c \in \mathbb{R}$, which we may assume nonnegative, we have, for all $n \in \mathbb{N}_0$,

$$0 \leq \alpha_n \leq \frac{c}{n^2} \quad \text{if } \alpha_n \geq 0,$$

$$0 \leq -\alpha_n \leq \frac{c}{n} \quad \text{if } \alpha_n \leq 0.$$

Therefore the sequence $\alpha_n$ tends to zero, independently of the sign assumed on each $I_n$. Since $X$ is unbounded, we are not allowed to buy or sell one unit of the risky investment $X$, and hence $X(1)$ is not a replicable claim.

We are now ready to analyze $M_1$. Every $Q \in M_1$ is identified by its density on $J_n^i$, denoted by $q_i(n)$. From the definition of $M_1$ in (3) we see that each $Q \in M_1$ is characterized by

$$\sum_{n \geq 1} \frac{q_1(n) + q_2(n)}{2^{n+1}} = 1 \quad \text{and} \quad q_1(n) = n q_2(n) \quad \forall n \geq 1,$$

which imply in particular that $\sum_{n \geq 1} \frac{(n+1)q_2(n)}{2^{n+1}}$ is finite. For later considerations, we observe also that $X(1)$ is not integrable for every $Q \in M_1$. Consider now the claim

$$f = \begin{cases} 1, & \text{on } J_n^1, \\ -n, & \text{on } J_n^2. \end{cases}$$

It is evident that $f \in L^1(Q)$ and $E_Q[f] = 0$ for any $Q \in M_1$. By using the duality relation in (10), we see that the weak super replication price of $f$ is equal to zero: $\hat{f}_{id} = 0$. However, $\hat{f} = 1$. Indeed if we try to write $f - x$ as $\alpha X(1) - h$ with $\alpha$ admissible and $h$ nonnegative, we obtain that, for every $n \geq 1$, the following must hold:

$$1 = n\alpha_n - h_1(n) + x,$$
$$-n = -n^2 \alpha_n - h_2(n) + x,$$

where $h_i(n)$ stands for the value of $h$ on $J_n^i$. Clearly the second equation can be always satisfied, provided that we choose $h_2(n)$ big enough.

Then analyzing the first one we get

$$h_1(n) = n\alpha_n + x - 1 \geq 0 \quad \forall n,$$

that is, $x \geq 1 - n\alpha_n$. Now, if $(\alpha_n)_n$ is definitely negative, we obviously get $x \geq 1$. In case $\alpha_n \geq 0$ infinitely many times, for these $\alpha_n$ we have $0 \leq \alpha_n \leq \frac{c}{n^2}$ and so $n\alpha_n$ is infinitesimal, when nonnegative. The consequence is again $x \geq 1$. Since $(f - 1) \in -L_+^0$, then $\hat{f} \leq 1$ and therefore $\hat{f} = 1$.

The difference between these two super replication prices is due to the fact that $f$ is equal to $(1, \frac{1}{2}, \frac{1}{3}, \ldots, \frac{1}{n}, \ldots)X(1)$, which is in $K^0 \cap \bigcap_{Q \in M_1} L^1(Q)$.

Under each $Q \in M_1$, this claim can be arbitrarily well $L^1(Q)$-approximated by claims in the form: $(1, \frac{1}{2}, \frac{1}{3}, \ldots, \frac{1}{n}, 0, 0, \ldots)X(1)$, which are in $K$ and have zero cost. When we require the usual stronger, pointwise condition $f - x = \alpha X(1) - h$, we obtain, due to the "artificial" admissibility requirement, the higher value $\hat{f} = 1$.

The difference between the weak and the classical super replication prices becomes more evident if we consider the claim $(kf)$ with $k \in \mathbb{R}$ positive and arbitrarily large. Reasoning exactly as before, we get $\widehat{(kf)} = k$. Selling at such an expensive price could be difficult, whereas the weak super replication price $\widehat{(kf)}_{id}$ is still zero. The drawback is that in this case one has to accept the possibility of only approximating $(kf - x)$ via bounded super-replicable claims in $C$.

EXAMPLE 9. Consider the same setup as in Example 8 and choose $\Phi(y) = y^2$, for $y \geq 0$. If we take $X(1)$ as the claim under consideration, it is rather easy to see that $\widehat{X(1)} = +\infty$, while $\sup_{Q \in M_1} E_Q[X(1)]$ is not even well defined.

In spite of these negative facts, the condition $E[\Phi(\frac{dQ}{dP})] < +\infty$ implies that $\sum_{n \geq 1} \frac{(n^2+1)q_2^2(n)}{2^{n+1}}$ is finite, thus $\{nq_2(n)2^{-(n+1)/2}\}_n \in l^2$. By the obvious remark $\{n2^{-(n+1)/2}\}_n \in l^2$, we get

$$\sum_{n \geq 1} \frac{n^2 q_2(n)}{2^{n+1}} < +\infty,$$

which, up to a constant, is just the $Q$-integrability condition on $X(1)$. Therefore, $X(1)$ is integrable for every $Q \in M_\Phi$ and the integral is zero. Summing up, we have

$$\widehat{X(1)}_\Phi = \sup_{Q \in M_\Phi} E_Q[X(1)] = 0 < \widehat{X(1)} = +\infty.$$

**2. Abstract formulation.** Recall that a subset $G$ of a vector space is a convex cone if $x, y \in G$ implies that $\alpha x + \beta y \in G$ for all $\alpha, \beta \geq 0$. Let $L \subseteq X, L' \subseteq X'$ be two convex cones in two vector spaces $X$ and $X'$. Let

$$\langle \cdot, \cdot \rangle : L \times L' \to \mathbb{R} \cup \{+\infty\}$$

be a "positive bilinear" form; that is, both applications $x \to \langle x, x' \rangle$ and $x' \to \langle x, x' \rangle$ are additive, positively homogeneous and equal to 0 at 0. We shall set $\langle x, x' \rangle \triangleq x'(x)$, for $x \in L$ and $x' \in L'$. With respect to $(L, L', \langle \cdot, \cdot \rangle)$ we define the polar $G^0$ and the bipolar $G^{00}$ of a convex cone $G$ by

$$G^0 \triangleq \{z \in L' | z(g) \leq 0 \; \forall g \in G\},$$
$$G^{00} \triangleq \{g \in L | z(g) \leq 0 \; \forall z \in G^0\}.$$

We assume that there exists an element, denoted by $\mathbf{1}$, such that $\mathbf{1} \in L$ and $-\mathbf{1} \in L$.

THEOREM 10. *Let $G \subseteq L$ be a convex cone satisfying $G^{00} = G$ and $-\mathbf{1} \in G$. If the set $N_1 \triangleq \{z \in G^0 | z(\mathbf{1}) = 1\}$ is not empty, then for all $f \in L$ we have*

$$\hat{f} \triangleq \inf\{x \in \mathbb{R} | f - x\mathbf{1} \in G\} = \sup\{z(f) | z \in N_1\}. \tag{13}$$

*In case $\hat{f} < +\infty$, it is a minimum.*

PROOF. First note that since $\mathbf{1} \in L$ and $-\mathbf{1} \in L$, then from $z(0) = 0$ and the additivity of all $z \in L'$ we deduce that $-\infty < z(-\mathbf{1}) = -z(\mathbf{1}) < +\infty$ and $z(f - x\mathbf{1})$ is well defined for all $z \in L'$, $f \in L$ and $x \in \mathbb{R}$. Hence $z(f - x\mathbf{1}) = z(f) - x$ for all $z \in N_1$ and $x \in \mathbb{R}$. Given $f \in L$ set $f^* \triangleq \sup\{z(f) | z \in N_1\} \leq +\infty$.

For all $x \in \mathbb{R}$ such that $(f - x\mathbf{1}) \in G$ we have $0 \geq \sup\{z(f - x\mathbf{1}) | z \in N_1\} = \sup\{z(f) | z \in N_1\} - x$ and hence $f^* \leq \hat{f}$.

To prove that $\hat{f} \leq f^*$ we may assume that $f^* < +\infty$ and it is sufficient to show that $(f - f^*\mathbf{1}) \in G$. Define

$$N \triangleq G^0 = \{z \in L' | z(g) \leq 0 \; \forall g \in G\} \tag{14}$$

and $N_0 \triangleq \{z \in N | z(\mathbf{1}) = 0\}$, so that $N = \bigcup_{\lambda > 0} \lambda N_1 \cup N_0$.

By definition of $f^*$, $-\infty < z(f - f^*\mathbf{1}) \leq 0$ for all $z \in N_1$. Let $z_0 \in N_0$ and note that if $z \in N_1$, then $(z + \lambda z_0) \in N_1$ for all $\lambda > 0$ and

$$0 \geq (z + \lambda z_0)(f - f^*\mathbf{1}) = z(f - f^*\mathbf{1}) + \lambda z_0(f) \quad \text{for all } \lambda > 0.$$

This implies $\lambda z_0(f) \leq -z(f - f^*\mathbf{1}) < +\infty$ for all $\lambda > 0$ and so $z_0(f) \leq 0$. Hence, $z_0(f - f^*\mathbf{1}) = z_0(f) \leq 0$ for all $z_0 \in N_0$. Therefore, $z(f - f^*\mathbf{1}) \leq 0$ for all $z \in N$ and we deduce that $(f - f^*\mathbf{1})$ belongs to the polar of $N$; that is, it belongs to $G^{00} = G$. □

REMARK 11. Note that the assumption that $N_1$ is not empty excludes that $\mathbf{1} = \mathbf{0}$. In our applications of Theorem 10, we will always consider $L \subseteq L^0$, $L' \subseteq L^1(P)$, $G$ will always be a convex cone containing $-L^\infty_+$, which implies that $N \triangleq G^0 \subseteq L^1_+$, and the element $\mathbf{1}$ will be the indicator function of $\Omega$. As a consequence of these conditions, $N_0 = \{0\}$.

REMARK 12. If $(L, L')$ is a dual system of vector spaces and if $\tau$ is any topology compatible with $(L, L')$, then the bipolar theorem, applied to the locally convex topological vector space $(L, \tau)$, guarantees $G^{00} = G$, whenever $G$ is a convex $\tau$-closed set.

2.1. *Proof of Theorem 3.*

DEFINITION 13 (see [4, 18]). A subset $C \subseteq L^0$ is Fatou closed if for every sequence $f_n \in C$ that is uniformly bounded from below and that converges $P$-a.s. to $f$, we have $f \in C$.

We collect in the following theorem some relevant results taken from Delbaen and Schachermayer (see [4, 5]).

THEOREM 14. (a) *If $D \subseteq L^0$ is a convex Fatou closed set, then $D \cap L^\infty$ is $\sigma(L^\infty, L^1)$-closed ([4], Theorem 4.2).*
(b) *If $M_1 \cap \mathbb{P} \neq \varnothing$, then $(K - L^0_+)$ is Fatou closed ([4], Theorem 4.2, and [5], Theorem 4.1).*

In [3] a bipolar theorem for $(L^0_+, L^0_+)$ is shown to hold, provided that the bilinear form $\langle \cdot, \cdot \rangle$ is allowed to take the value $+\infty$. The proof of Theorem 15(a) is based on the proof of the simpler bipolar theorem for $(L^{bb}, L^1_+)$ in [12].

THEOREM 15. (a) *If $C_{bb}$ is Fatou closed, then $C_{bb} = (C_{bb})^{00}$.*
(b) *In particular if $M_1 \cap \mathbb{P} \neq \varnothing$, then $C_{bb} = (C_{bb})^{00}$.*

PROOF. By definition, $(C_{bb})^0 \triangleq \{z \in L^1_+ : E[zf] \leq 0 \,\forall f \in C_{bb}\}$ and $(C_{bb})^{00} \triangleq \{f \in L^{bb} : E[zf] \leq 0 \,\forall z \in (C_{bb})^0\}$.

(a) Clearly $C_{bb} \subseteq (C_{bb})^{00}$. To show that $(C_{bb})^{00} \subseteq C_{bb}$ suppose by contradiction that there exists $f \in (C_{bb})^{00}$ and $f \notin C_{bb}$. Then $f_n \triangleq (f \wedge n) \in (C_{bb})^{00} \cap L^\infty$, $f_n \uparrow f$ $P$-a.s. and $f_n$ is uniformly bounded from below. Since $C_{bb}$ is Fatou closed and $f \notin C_{bb}$, then there exists $n_0$ such that $f_{n_0} \notin C_{bb}$. Since the set $C_{bb} \cap L^\infty$ is convex and $\sigma(L^\infty, L^1)$-closed [see Theorem 14(a)] and $f_{n_0} \notin C_{bb} \cap L^\infty$ the separation theorem in $(L^\infty, \sigma(L^\infty, L^1))$ guarantees the existence of $z \in L^1$ such that

$$E[zg] \leq 0 \quad \forall g \in C_{bb} \cap L^\infty \quad \text{and} \quad E[zf_{n_0}] > 0.$$

Since $-L^\infty_+ \subseteq C_{bb} \cap L^\infty$ we have $z \in L^1_+$. We now show that $z \in (C_{bb})^0$, which is in contradiction with $f_{n_0} \in (C_{bb})^{00}$ and $E[zf_{n_0}] > 0$. For each $\tilde{g} \in C_{bb}$ we set $g_n \triangleq (\tilde{g} \wedge n)$. Then $g_n \in C_{bb} \cap L^\infty$, $g_n \uparrow \tilde{g}$, $P$-a.s. and $g_n$ is uniformly bounded from below. By Fatou's lemma,

$$E[z\tilde{g}] \leq \lim E[zg_n] \leq 0 \quad \forall \tilde{g} \in C_{bb}.$$

(b) From Theorem 14(b) we know that $(K - L^0_+)$ is Fatou closed; hence $C_{bb} = (K - L^0_+) \cap L^{bb}$ is Fatou closed and (b) follows from (a). □

Now we are ready to give a proof, based on Theorem 10, of Theorem 3.

PROOF OF THEOREM 3. To prove (9), we apply Theorem 10, with $L = L^{bb}$, $L' = L^1_+$, $\mathbf{1} = \mathbb{1}_\Omega$ and $G = C_{bb}$. The positive bilinear form will be $x'(x) = E[x'x]$.

From (14) we get $N \triangleq (C_{bb})^0 = \{z \in L^1_+ | E[zg] \leq 0 \ \forall g \in C_{bb}\}$ and $N_1 \triangleq \{z \in N | E[z] = 1\}$. Since

$$\{z \in L^1_+ | E[zg] \leq 0 \ \forall g \in C_{bb}\} = \{z \in L^1_+ | E[zg] \leq 0 \ \forall g \in K\},$$

we may identify $N_1$ with $M_1$. From Theorem 15(b) we see that the assumptions of Theorem 10 are satisfied. Hence

$$\inf\{x \in \mathbb{R} | f - x \in C_{bb}\} = \sup\{E[zf] | z \in M_1\}. \qquad \square$$

**3. The polarity between $C_\Phi$ and $co(M_\Phi)$.** In this section we stick to the terminology of [11], Chapter 8. Define the linear spaces

$$L = \bigcap_{Q \in M_\Phi} L^1(Q) \quad \text{and} \quad L' = Lin\{M_\Phi\} \subseteq L^1(P),$$

where we assume that $M_\Phi$ is not empty and we identify each $Q$ with its Radon–Nikodym derivative w.r.t. $P$.

Notice that $C \subseteq L^\infty(P) \subseteq L$. For all $z \in L$ and $z' \in L'$, we have that $(zz') \in L^1(P)$ and the bilinear form $z \times z' \to E[zz']$ is well defined. Then $(L, L')$ defines a dual system.

DEFINITION 16. We denote by $\tau$ a locally convex topology on $L$ compatible with the duality $(L, L')$.

Just by definition, endowed with the $\tau$-topology $L$ is a locally convex topological vector space where the set of continuous linear forms on $L$ is precisely $L'$. We may select any topology compatible with the dual system $(L, L')$, since our results depend only on the property that the topological dual of $L$ is $L'$.

Note that this topology $\tau$ needs not to be Hausdorff, since generally $L'$ does not separate points in $L$. Think of the case when we have just one element in $M_\Phi$ (a complete market case, in which the unique equivalent pricing measure has finite entropy).

Define

$$(15) \qquad C_\Phi \triangleq \bigcap_{Q \in M_\Phi} \overline{(K - L^1_+(Q))}^Q.$$

The main result of this section is the following theorem. Its proof will be based on Proposition 19 and Theorem 20, which will also provide a different representation for $C_\Phi$.

THEOREM 17.   *Assume that $\Phi(0) < +\infty$ and $M_\Phi \cap \mathbb{P} \neq \varnothing$. With respect to the topology $\tau$ we have:* (a) *$C_\Phi$ is the closure of $C$;* (b) *$C_\Phi$ and the convex cone $co(M_\Phi)$ generated by $M_\Phi$ are polar to one another.*

As an immediate consequence of Theorems 10 and 17 we prove Theorem 5.

PROOF OF THEOREM 5.   Since $M_\Phi \subseteq M_1$, the inequality in (11) is proved in Proposition 1. Consider the dual system $(L, L')$ and the topology $\tau$ on $L$. Set $G = C_\Phi$. From Theorem 17 we deduce $N = (C_\Phi)^0 = co(M_\Phi)$ and $N_1 = M_\Phi$. The assumptions of Theorem 10 are satisfied and then from (13) we get

$$\inf\{x \in \mathbb{R} | f - x \in C_\Phi\} = \sup\{E_Q[f] | Q \in M_\Phi\}. \qquad \square$$

PROPOSITION 18.   *Assume that $\Phi(0) < +\infty$. If $Q_0 \ll P$, $Q_1 \ll P$, $x \in (0,1)$, $Q = xQ_1 + (1-x)Q_0$, then*

$$E\Phi\left(\frac{dQ}{dP}\right) < +\infty \quad \text{if and only if}$$

$$E\Phi\left(\frac{dQ_0}{dP}\right) < +\infty \quad \text{and} \quad E\Phi\left(\frac{dQ_1}{dP}\right) < +\infty.$$

PROOF.   The convexity of $\Phi$ implies that $E\Phi(\frac{dQ}{dP}) < +\infty$ if $E\Phi(\frac{dQ_i}{dP}) < +\infty$ for $i = 0, 1$. Conversely suppose that $E\Phi(\frac{dQ}{dP}) < +\infty$. For $i = 0, 1$, we have $\Phi^-(\frac{dQ_i}{dP}) \in L^1(P)$, since $\Phi$ is convex and $\frac{dQ_i}{dP} \in L^1(P)$. Therefore we only need to show the integrability of $\Phi^+(\frac{dQ_i}{dP})$, which is trivially true if $\Phi(+\infty) < +\infty$. If $\Phi(+\infty) = +\infty$ then $\Phi^+$ is nondecreasing on $(y_0, +\infty)$ for some $y_0 > 0$. From $Q = xQ_1 + (1-x)Q_0$ we deduce

$$\frac{dQ_1}{dP} = \frac{1}{x}\frac{dQ}{dP} - \frac{1-x}{x}\frac{dQ_0}{dP} \leq \frac{1}{x}\frac{dQ}{dP}, \qquad P\text{-a.s.},$$

$$E\Phi^+\left(\frac{dQ_1}{dP}\right) = E\left[\Phi^+\left(\frac{dQ_1}{dP}\right)\mathbb{1}_{\{dQ_1/dP \leq y_0\}}\right] + E\left[\Phi^+\left(\frac{dQ_1}{dP}\right)\mathbb{1}_{\{dQ_1/dP > y_0\}}\right]$$

$$\leq \max_{0 \leq y \leq y_0} \Phi^+(y) + E\left[\Phi^+\left(\frac{1}{x}\frac{dQ}{dP}\right)\mathbb{1}_{\{dQ_1/dP > y_0\}}\right] < +\infty$$

since, from the growth condition $G(\Phi)$, we have $\Phi^+(\frac{1}{x}\frac{dQ}{dP}) \leq \alpha\Phi^+(\frac{dQ}{dP}) + \beta(\frac{dQ}{dP} + 1) \in L^1(P)$. Similarly for $\frac{dQ_0}{dP}$.   $\square$

Let $\overline{C}$ be the closure of $C$ with respect to the $\tau$ topology. Note that $\overline{C}$ is a convex cone and $\overline{C} \subseteq L \subseteq L^1(Q)$ for all $Q \in M_\Phi$. The polar of $\overline{C}$ with respect to the $\tau$ topology is given by

$$\overline{C}^0 \triangleq \{z' \in L' : E[zz'] \leq 0 \text{ for all } z \in \overline{C}\} \subseteq L^1_+(P),$$

since $-L_+^\infty \subseteq C$.

PROPOSITION 19. *If $\Phi(0) < +\infty$, then $co\{M_\Phi\} = \overline{C}^0$.*

PROOF. All $Q \in M_\Phi$ are $\tau$-continuous linear functionals, so that (for a fixed $Q$) the set $\{z \in L | E_Q[z] \leq 0\}$ is $\tau$-closed and it contains $C$. We deduce that if $z \in \overline{C}$, then $E_Q[z] \leq 0$ for all $Q \in M_\Phi$. Since $M_\Phi$ is convex, $L'$ admits the following representation:

$$L' = \{z' \in L^1(P) : z' = \alpha z_1' - \beta z_0', \ \alpha, \beta \geq 0, \ z_1', z_0' \in M_\Phi\}.$$

We claim that $M_\Phi = \overline{C}_1^0 \triangleq \overline{C}^0 \cap \{\text{unit sphere of } L^1(P)\}$. Note that

$$\overline{C}_1^0 = \{Q \ll P : Q = (1+\beta)Q_1 - \beta Q_0, \ \beta \geq 0, \ Q_1, Q_0 \in M_\Phi$$
$$\text{and } \forall z \in \overline{C}, \ E_Q[z] \leq 0\}.$$

Obviously $M_\Phi \subseteq \overline{C}_1^0$: so we consider the case $\beta > 0$. If $Q \in \overline{C}_1^0$, then $\forall z \in \overline{C}$, $E_Q[z] \leq 0$ and so $Q \in M_1$. It remains only to check that if $Q \in \overline{C}_1^0$, then $E\Phi(\frac{dQ}{dP}) < +\infty$. If $Q \triangleq (1+\beta)Q_1 - \beta Q_0$, then $Q_1 = \frac{1}{1+\beta}Q + \frac{\beta}{1+\beta}Q_0 = xQ + (1-x)Q_0$, $x = \frac{1}{1+\beta} \in (0,1)$, and the thesis follows from Proposition 18. □

The following theorem is proved in [9], Theorem 3 adding to $G(\Phi)$ the assumptions that $\Phi$ is strictly convex and differentiable. But the proof of the theorem remains unchanged even without these additional assumptions. Let

$$(co(M_\Phi))^0 \triangleq \{f \in L : E_Q[f] \leq 0 \ \forall Q \in M_\Phi\}.$$

THEOREM 20. *If $M_\Phi \cap \mathbb{P} \neq \varnothing$, then*

$$C_\Phi = \bigcap_{Q \in M_\Phi} \overline{C}^Q = (co(M_\Phi)^0).$$

PROOF OF THEOREM 17. Since $co\{M_\Phi\} = \overline{C}^0$, the bipolar $\overline{C}^{00}$ of $\overline{C}$ is given by:

$$\overline{C}^{00} \triangleq \{z \in L : E[zz'] \leq 0 \text{ for all } z' \in \overline{C}^0\}$$
$$= \{z \in L : E_Q[z] \leq 0 \text{ for all } Q \in M_\Phi\} = C_\Phi,$$

by Theorem 20. From the bipolar theorem we deduce that $\overline{C} = \overline{C}^{00} = C_\Phi$. From $co\{M_\Phi\} = \overline{C}^0$ we then get $(co\{M_\Phi\})^0 = C_\Phi$ and $(C_\Phi)^0 = co\{M_\Phi\}$. □

The boundedness of $\Phi$ in a right neighborhood of 0 is essential in Propositions 18 and 19 and in Theorem 17, as the following example shows.

EXAMPLE 21. The context is the same of Example 8. Consider the function $\Phi$ defined by:

$$\Phi = \begin{cases} -\ln(y), & \text{on } 0 < y \leq 1, \\ y^2 - 3y + 2, & \text{on } y > 1. \end{cases}$$

Obviously, $\Phi$ is strictly convex and differentiable. The point is that in this model there exists a $Q_1 \in M_1$, with $Q_1$ not equivalent to $P$ and with bounded density: such a measure has infinite generalized entropy, that is, $Q_1 \notin M_\Phi$. For instance, let $\frac{dQ_1}{dP} = 2\chi_{I_1} = 2\chi_{(\frac{1}{2},1]}$. Then, pick any $Q_0 \in M_\Phi$: for example, take $\frac{dQ_0}{dP}$ equal to $c\frac{n}{e^n}$ on $J_n^1$ (and consequently equal to $\frac{c}{e^n}$ on $J_n^2$), where $c$ is the normalizing constant. Consider now the convex combination $Q^x = (1-x)Q_0 + xQ_1$, $x \in (0,1)$. Since the following inequalities hold true

$$(1-x)Q_0 \leq Q^x \leq (1-x)Q_0 + \text{const},$$

$Q^x$ has finite generalized entropy, that is, $Q^x \in M_\Phi$.

Since $Q_1 \notin M_\Phi$, to show that $co(M_\Phi) \subsetneq (C_\Phi)^0$ it is sufficient to show that $Q_1 \in (C_\Phi)^0$. It is obvious that $Q_1 \in Lin(M_\Phi) = L'$ and $\overline{C} \subseteq L^1(Q_1)$. Recall that $\overline{C} = C_\Phi \subseteq \overline{C}^Q$ and $E_Q[f] \leq 0$ for all $Q \in M_\Phi$ and $f \in C_\Phi$. Since $|f|\frac{dQ^x}{dP} \leq |f|(\frac{dQ_0}{dP} + \frac{dQ_1}{dP})$ we deduce, if $f \in C_\Phi$, $E_{Q_1}[f] = \lim_{x \to 1} E_{Q^x}[f] \leq 0$.

REMARK 22. Motivated by the last lines of the previous example, we now make some extra observations on the duality $(L, L')$. As we have already noted, the dual system may not be separated. The consequence is that in general we cannot put a topology $\mu$ on $L'$ which is compatible with the duality $(L, L')$, that is, such that the dual of $(L', \mu)$ is exactly $L$ (think again of the case when $|M_\Phi| = 1$).

However, if we define on $L$ the equivalence relation $\sim$,

$$f \sim g \quad \text{iff } E_Q[f] = E_Q[g] \text{ for all } Q \in M_\Phi,$$

and we define $\frac{L}{\sim}$ to be the quotient of $L$ w.r.t. the relation $\sim$, then it can be easily seen that $\frac{L}{\sim}$ is a vector space with the obviously defined sum and scalar multiplication.

We indicate with $\tau_\sim$ the quotient topology of $(L, \tau)$ on $\frac{L}{\sim}$. It is now a simple exercise proving that, for all $\xi \in \frac{L}{\sim}$ and $z' \in L'$, we have that $zz' \in L^1(P)$ (where $z$ is a generic element of the equivalence class $\xi$) and the bilinear form $\xi \times z' \to \prec \xi, z' \succ \triangleq E[zz']$ is well defined. Then $(\frac{L}{\sim}, L')$ is a dual system, it is separating and the topology $\tau_\sim$ on $\frac{L}{\sim}$ is compatible. Now we also can endow $L'$ with a topology $\nu$ compatible with this new system.

When the condition $\Phi(0) < +\infty$ is satisfied, we have that $co(M_\Phi)$ coincides with $(\frac{C_\Phi}{\sim})^0$ and therefore is $\nu$-closed.

The previous example shows that this is not always the case when $\Phi(0)$ is infinite. In fact, fix an $\eta \in \frac{L}{\sim}$. Then, with the same notation used before, $\prec \eta, Q^x \succ$ tends to $\prec \eta, Q^1 \succ$ when $x \to 1$. Now, letting $\eta$ vary arbitrarily in $\frac{L}{\sim}$ we get that $Q^x$ tends to $Q^1$ in the $\nu$-topology. Therefore neither $M_\Phi$ nor $co(M_\Phi)$ is $\nu$-closed.

**4. Comparison with the Delbaen–Schachermayer approach, when $\Phi = id$.** In their remarkable paper [5], Delbaen and Schachermayer introduced the notions of feasible weight function $w$ for the process $X$ and of $w$-admissible integrands for $X$ to get the duality results stated below in Theorem 25. We recall here some of their definitions and results and we defer to [5], Section 5, for their motivation and explanation. In the sequel it is always assumed that $M_1 \cap \mathbb{P} \neq \varnothing$. Note also that the time horizon $T$ appearing throughout this paper could be finite as well as $+\infty$: the latter case will be now considered.

DEFINITION 23 ([5], Definition 5.1). If $w \geq 1$ is a random variable, if there is $Q_0 \in M_\sigma \cap \mathbb{P}$ such that $E_{Q_0}[w] < \infty$, then we say that the integrand $H$ is $w$-admissible if there exists some nonnegative real number $c$ such that, for each element $Q \in M_\sigma \cap \mathbb{P}$ and each $t \geq 0$, we have that $(H \cdot X)_t \geq -cE_Q[w|\mathcal{F}_t]$.

DEFINITION 24 ([5], Definition 5.4). A real random variable $w \geq 1$ is called a feasible weight function for $X$ if the following hold:

(a) there is a strictly positive bounded predictable process $\phi$ such that the maximal function of the $\mathbb{R}^d$-valued stochastic integral $\phi \cdot X$ satisfies $(\phi \cdot X)^* \leq w$;

(b) there is an element $Q_0 \in M_\sigma \cap \mathbb{P}$ such that $E_{Q_0}[w] < \infty$.

As pointed out in the cited article, feasible weight functions do exist. Let $w$ be a feasible weight function for $X$ and set

$$K_w \triangleq \{(H \cdot X)_\infty | H \text{ is } w\text{-admissible}\},$$

$$\hat{f}_w \triangleq \inf\{x \in \mathbb{R} | f - x \in K_w - L^0_+\},$$

$$M_{\sigma,w} \triangleq \{Q \in M_\sigma | E_Q[w] < \infty\}.$$

THEOREM 25 ([5], Theorem 5.5). *If $w$ is a feasible weight function and $f$ is a random variable such that $f \geq -w$, then*

(16) $$\hat{f}_w = \inf\{x \in \mathbb{R} | f - x \in K_w - L^0_+\} = \sup_{Q \in M_{\sigma,w} \cap \mathbb{P}} E_Q[f]$$

*and if the quantities are finite, the infimum is a minimum.*

We now compare the super replication price $\hat{f}_w$ of $f$ given in (16) with the weak super replication price $\hat{f}_{id}$ of $f$ given in (10).

The first important remark is that given a claim $f \in \bigcap_{Q \in M_1} L^1(Q)$ then $\hat{f}_{id}$ is uniquely defined and is not dependent on the agent. On the contrary, the super replication price $\hat{f}_w$, of the same claim $f$, *will in general depend on the different feasible weight functions $w$* selected by the investor. Indeed, $\hat{f}_w$ depends on how much one is ready to lose in the trading. By admitting bigger losses, this price decreases, as we will show in the example in Section 4.1. Only admitting the knowledge of a feasible weight function $w$, the super replication price $\hat{f}_w$ of those claims $f$ satisfying $f \geq -w$ is uniquely defined and (16) may be applied.

If $f \in \bigcap_{Q \in M_1} L^1(Q)$, then by simply considering $w(f) \triangleq w \vee f^-$ (where $f^-$ is the negative part of $f$) we obtain a feasible weight function such that $f \geq -w(f)$. Therefore, for each given claim $f \in \bigcap_{Q \in M_1} L^1(Q)$ we can always find at least one suitable feasible weight $w_f$ so that we can apply the duality formula (16) to the couple $f, w_f$ to get the particular super replication price $\hat{f}_{w_f}$.

From (16), (10) and Remark 2, we get

$$f \in \bigcap_{Q \in M_1} L^1(Q) \quad \implies \quad \hat{f}_{id} = \sup_{Q \in M_1 \cap \mathbb{P}} E_Q[f] \geq \sup_{Q \in M_{\sigma,w_f} \cap \mathbb{P}} E_Q[f] = \hat{f}_{w_f}.$$

In [5] it is also proved that $M_\sigma \cap \mathbb{P}$ is dense in $M_1 \cap \mathbb{P}$ (Proposition 4.7) and that $M_{\sigma,w} \cap \mathbb{P}$ is dense in $M_\sigma \cap \mathbb{P}$ (Corollary 5.13). Unfortunately, in spite of the density properties, we cannot apply the dominated convergence theorem, as done in Remark 2. As shown in Example 29, the weak super replication price $\hat{f}_{id}$ can be strictly greater than $\hat{f}_{w(f)}$ (or than $\hat{f}_w$ with any $w$ feasible with $f \geq -w$).

4.1. *Dependence on $w$*. First recall that for locally bounded processes, as those we will consider in this section, the sets $M_1$ of separating measures and $M_\sigma$ of $\sigma$-martingale measures are equal and coincide with the set of local martingale measures. Hence $M_{1,w} \triangleq \{Q \in M_1 | E_Q[w] < \infty\} = M_{\sigma,w}$ and $M_1$ may replace $M_\sigma$ (and vice versa) in any subsequent formulas.

With the next example we provide evidence of the dependence of the super replication price $\hat{f}_w$ from the feasible weight function $w$ and of a situation in which

$$(17) \qquad \sup_{Q \in M_\sigma \cap \mathbb{P}} E_Q[f] > \sup_{Q \in M_{\sigma,w} \cap \mathbb{P}} E_Q[f].$$

Example 5.14 in [5] was exactly intended to prove the previous inequality, but, as we now explain, it is not correct. The claim $f$ and the feasible weight function $w_1$, introduced in the next example, are exactly those considered in Example 5.14 in [5]. However, we will prove in item 5 below [see also (23)] that, contrary to the assertion (2) made after Example 5.14 in [5], the two suprema in (17) coincide for such $f$ and $w_1$. For the validity of the strict inequality in (17) (or in [5], (5.1)) we have to use a different weight function ($w_2$) and to exploit the peculiar feature (see Lemma 27) of a positive strict local martingale $X$ under $P$, which admits a probability measure $Q \sim P$ such that $X \in H^2(Q)$.

EXAMPLE 26. On a suitable stochastic basis $(\Omega, (\mathcal{F}_t)_{t\geq 0}, P)$ there exist:

(a) a continuous process $S$ satisfying $S_0 = 0$ such that $P \in M_1 \cap \mathbb{P}$, where $M_1$ is the set of separating measures for $S$;
(b) two $S$-feasible weight functions $w_1$ and $w_2$;
(c) a claim $f \in \bigcap_{Q \in M_1} L^1(Q)$ satisfying $f \geq -w_1$, $f \geq -w_2$;

such that:

1. $w_1 \in \bigcap_{Q \in M_1} L^1(Q)$, so that $M_{\sigma, w_1} = M_\sigma = M_1$;
2. $S$ is uniformly bounded from above and is a submartingale for each $Q \in M_1$;
3. $S$ is not a martingale under $P$ and $E_P[S_\infty] > 0$;
4. $\forall R \in M_{\sigma, w_2}$, $S$ is an $R$-uniformly integrable martingale and $E_R[S_\infty] = 0$;
5. $\hat{f}_{id} = \hat{f}_{w_1} > \hat{f}_{w_2} = 0$.

To demonstrate this example, we need a result based on a slight modification of the example in [6], Section 2, to which we refer for a detailed construction.

We call
$$L_t \triangleq \exp(B_t - \tfrac{1}{2}t)$$
and
$$(18) \qquad N_t^{(a)} \triangleq \exp\left(aW_t - \frac{a^2}{2}t\right),$$

where $a$ is a positive real constant and $(B, W)$ is a standard two-dimensional Brownian motion on a stochastic basis $(\Omega, (\mathcal{F}_t)_{0 \leq t \leq +\infty}, P)$. We assume that the filtration $\mathcal{F}$ is the augmentation of the natural one, $(\mathcal{F}_t^{B,W})_t$, induced by $(B, W)$. Both $L$ and $N^{(a)}$ are positive, strict $P$-local martingales. Then, define the stopping times

$$(19) \qquad \tau \triangleq \inf\{t | L_t = \tfrac{1}{2}\},$$

$$(20) \qquad \sigma^{(a)} \triangleq \inf\{t | N_t^{(a)} = 2\}.$$

Notice that
$$\tau = \inf\{t | B_t - \tfrac{1}{2}t = \log \tfrac{1}{2}\},$$
$$\sigma^{(a)} = \inf\left\{t | W_t - \frac{a}{2}t = \frac{\log 2}{a}\right\},$$

so these two stopping times are passage times of Brownian motion with drift.

Now define the stopped processes $X^{(a)} \triangleq L^{\tau \wedge \sigma^{(a)}}$ and $Y^{(a)} \triangleq (N^{(a)})^{\tau \wedge \sigma^{(a)}}$ and the probability measure $Q^{(a)} \triangleq Y^{(a)}_\infty \cdot P$.

The following result is analogous to Theorem 2.1 of [6], but the introduction of the parameter $a$ in (18) allows us to add item (d). When $a = 1$, Lemma 27 reduces to Theorem 2.1 of [6]. However, $X^{(1)}$ is not in $H^2(Q^{(1)})$.

LEMMA 27. (a) *For every $a > 0$, the process $X^{(a)}$ is a continuous strict local martingale under $P$ and $X^{(a)}_\infty > 0$ a.s., $X^{(a)}_0 = 1$, $E_P[X^{(a)}_\infty] < 1$.*

(b) *For every $a > 0$, the process $Y^{(a)}$ is a continuous uniformly bounded integrable martingale, that is strictly positive on $[0, +\infty]$.*

(c) *For every $a > 0$, the process $X^{(a)}$ is a uniformly integrable martingale under $Q^{(a)}$.*

(d) *$X^{(a)}$ belongs to $H^2(Q^{(a)})$ iff $a^2 \geq 8$.*

PROOF. We only need to prove item (d) since the first three points can be easily checked as in Theorem 2.1 of [6]. For simplicity of notation the dependence on $a$ is dropped.

By definition, $X$ is in $H^2(Q)$ iff $E_Q[\langle X \rangle_\infty] < +\infty$. Taking into account the positivity of the processes, an application of Doob's optional sampling theorem to the $P$-uniformly integrable martingale $N^\sigma$ leads to

$$E_Q[\langle X \rangle_\infty] = E[Y_\infty \langle L \rangle_{\tau \wedge \sigma}] = E[N_\sigma \langle L \rangle_{\tau \wedge \sigma}]$$

and, thanks to the independence of $(L, \tau)$ and $\sigma$, the last term becomes

$$(21) \qquad 2 \int \chi_{\{\sigma < +\infty\}}(\omega') E[\langle L \rangle_{\tau \wedge \sigma(\omega')}] \, dP(\omega').$$

Let us then analyze $E[\langle L \rangle_{\tau \wedge t}]$: it is equal to $E[L^2_{\tau \wedge t}]$ because $L^t$ is a square integrable martingale. By the Girsanov theorem we can write

$$E[\langle L \rangle_{\tau \wedge t}] = E[L^2_{\tau \wedge t}] = E[\exp\{2B_{\tau \wedge t} - \tau \wedge t\}] = \overline{E}[\exp\{\tau \wedge t\}],$$

where the last expectation is taken under the unique probability $\overline{P}$ on $\mathcal{F}^{B,W}_\infty$ such that $(\overline{B}_r)_r = (B_r - 2r)_r$ is a standard Brownian motion. With such a change of measure, $\tau = \inf\{r | \overline{B}_r + \frac{3}{2}r = -\log 2\}$ and the law of $\tau$ on $(0, +\infty]$ under $\overline{P}$ is given by

$$\tau(\overline{P}) = \frac{|b|}{\sqrt{(2\pi t^3)}} \exp\left[-\frac{(b - \mu t)^2}{2t}\right] dt + (1 - \exp(\mu b - |\mu b|))\varepsilon_{\{+\infty\}},$$

where $\mu = \frac{3}{2}$, $b = -\log 2$; that is, it consists of the sum of two positive measures, the first a.c. with respect to the Lebesgue measure on $(0, +\infty)$ with density

$$(22) \qquad f(t) = \frac{|b|}{\sqrt{(2\pi t^3)}} \exp\left[-\frac{(b - \mu t)^2}{2t}\right]$$

and the second being an atom in $+\infty$ with mass $1 - \exp(\mu b - |\mu b|)$ (see [16], page 196). Then

$$e^t \geq \overline{E}[\exp\{\tau \wedge t\}] = \int e^{s \wedge t} f(s)\, ds + \tfrac{7}{8} e^t \geq \tfrac{7}{8} e^t,$$

and the quantity in (21) is finite if and only if $E[\chi_{\{\sigma < +\infty\}} e^\sigma] < \infty$. Using the density $f(t)$ in (22), with $\mu = -\tfrac{a}{2}$, $b = \tfrac{\log 2}{a}$, of the absolutely continuous part of the law of $\sigma$ under $P$, we get

$$E[\chi_{\{\sigma < +\infty\}} e^\sigma] = \int_0^{+\infty} e^t \frac{(\log 2)/a}{\sqrt{(2\pi t^3)}} \exp\left[-\frac{((\log 2)/a + (a/2)t)^2}{2t}\right] dt$$

and the integral is finite iff $a^2 \geq 8$. $\square$

REMARK 28. Similar results can be obtained by replacing the constant $\tfrac{1}{2}$ in (19) with any $0 < c_1 < 1$ and the constant 2 in (20) with any $c_2 > 1$.

EXAMPLE 29 (Continued). Fix any $a > 0$ and take $X \triangleq X^{(a)}$, $P$, $Q \triangleq Q^{(a)}$ as defined before Lemma 27.

We define $S = 1 - X$. Then $P \in M_1$. We note that $S_0 = 0$ and $S$ is bounded from above, so that $H = -1$ is a "usual" admissible integrand. Under each $R \in M_1$, $-S$ is a supermartingale and hence $S$ is a submartingale.

We take $f = S_\infty$ as the claim to be evaluated. We are in a continuous context, so a $w \geq 1$ is feasible as soon as there exists a measure $R \in M_1 \cap \mathbb{P}$ such that $E_R[w]$ is finite.

First we consider $w_1 = 1 + X_\infty$. Note that $f \geq -w_1$ and that $w_1$ is feasible, since it is integrable for all $R \in M_1$ by construction. Note that when $a = 1$ this setting is precisely the one considered in Example 5.14 of [5]. Then the duality formula (16) can be applied to $f$ and we have, recalling Remark 2,

$$(23) \qquad \hat{f}_{id} = \sup_{R \in M_1 \cap \mathbb{P}} E_R[f] = \sup_{R \in M_{\sigma, w_1} \cap \mathbb{P}} E_R[f] = \hat{f}_{w_1} \geq E_P[f] > 0.$$

As a consequence of the last inequality, $H = 1$ is NOT $w_1$-admissible. If it were $S = (1 \cdot S)$ would become a supermartingale (this implication derives from Proposition 3.3 in [1] as well as from Theorem 5.3 in [5]) under each $R \in M_1$ and hence a martingale: this would imply $E_P[f] = 0$. Another argument is that, using the duality in (16), $\hat{f}(w_1) \leq 0$, a contradiction.

We now consider $w_2 = (X_\infty^*)^2$, where $X_t^* = \sup\{|X_s| | 0 \leq s \leq t\} = \sup\{X_s | 0 \leq s \leq t\}$. Now we need to assume that $a \geq 2\sqrt{2}$.

Then $w_2$ is certainly $Q$-integrable [by the Burkholder–Davis–Gundy inequalities, $w_2 \in L^1(Q)$; it is not in $L^1(P)$, because otherwise $X$ would be a $P$-square integrable martingale]: so, $w_2$ also is feasible and clearly $f \geq -w_2$. Now we get

$$\hat{f}_{w_2} = \sup_{R \in M_{1, w_2} \cap \mathbb{P}} E_R[f] = 0$$

because under these $R$ we obviously have

$$S_t = 1 - X_t \geq -E_R[w_2|\mathcal{F}_t];$$

that is, $H = 1$ is $w_2$-admissible and henceforth $S$ is an $R$-martingale. The crucial point that $M_{1,w_2} \cap \mathbb{P} \neq \varnothing$ was shown in Lemma 27, item 4.

**Acknowledgments.** We warmly thank the referees, an Associate Editor, F. Delbaen, W. Schachermayer and M. Schweizer for their suggestions, which helped us improve the paper.

Dipartimento di Economia  
Sezione di Finanza Matematica  
Università degli Studi di Perugia  
via A. Pascoli 20  
06123 Perugia  
Italy  
e-mail: s.biagini@unipg.it  

Dipartimento di Matematica  
per le Decisioni  
Università degli Studi di Firenze  
via C. Lombroso 6/17  
50134 Firenze  
Italy  
e-mail: marco.frittelli@dmd.unifi.it